\documentclass[12pt]{amsart}

\usepackage{amsthm}
\usepackage{amsmath}
\usepackage{amssymb}
\usepackage{eucal}
\usepackage{amsfonts}
\usepackage{amsmath}
\usepackage[Lenny]{fncychap}
\usepackage{enumerate}
\usepackage{amssymb}
\usepackage[english]{babel}
\usepackage[latin1]{inputenc}
\usepackage[breaklinks=true]{hyperref}
\usepackage[all]{xy}
\usepackage{fancybox}
\usepackage{latexsym,mathrsfs,longtable,lscape,dsfont,yfonts}
\usepackage{srcltx}

\ifx\pdfpagewidth\undefined 
\usepackage[T1]{fontenc}
\else 
\usepackage[OT1]{fontenc} 
\fi 
\usepackage{amsfonts,amssymb,latexsym,longtable,amsmath}

\newtheorem{thm}{Theorem}[section]
\newcommand{\bthm}{\begin{thm}} \newcommand{\ethm}{\end{thm}}
\newtheorem{prop}[thm]{Proposition}
\newcommand{\bprp}{\begin{prop}} \newcommand{\eprp}{\end{prop}}
\newtheorem{fact}[thm]{Fact}
\newcommand{\bfct}{\begin{fact}} \newcommand{\efct}{\end{fact}}
\newtheorem{prob}[thm]{Problem}
\newcommand{\bprb}{\begin{prob}} \newcommand{\eprb}{\end{prob}}
\newtheorem{quest}[thm]{Question}
\newcommand{\bqtn}{\begin{quest}} \newcommand{\eqtn}{\end{quest}}
\newtheorem{lem}[thm]{Lemma}
\newcommand{\blem}{\begin{lem}} \newcommand{\elem}{\end{lem}}
\newtheorem{claim}[thm]{Claim}
\newcommand{\bclm}{\begin{claim}} \newcommand{\eclm}{\end{claim}}
\newtheorem{cor}[thm]{Corollary}
\newcommand{\bcor}{\begin{cor}} \newcommand{\ecor}{\end{cor}}
\newtheorem{conj}[thm]{Conjecture}
\newcommand{\bcnj}{\begin{conj}} \newcommand{\ecnj}{\end{conj}}
\theoremstyle{definition}
\newtheorem{defn}[thm]{Definition}
\newcommand{\bdfn}{\begin{defn}} \newcommand{\edfn}{\end{defn}}
\newtheorem{spec}[thm]{Specializing}
\newcommand{\bspc}{\begin{spec}} \newcommand{\espc}{\end{spec}}
\theoremstyle{remark}
\newtheorem{rem}[thm]{Remark}
\newcommand{\brem}{\begin{rem}} \newcommand{\erem}{\end{rem}}
\newtheorem{cnv}[thm]{Convention}
\newcommand{\bcnv}{\begin{cnv}} \newcommand{\ecnv}{\end{cnv}}
\newtheorem{exam}[thm]{Example}
\newcommand{\bexm}{\begin{exam}} \newcommand{\eexm}{\end{exam}}
\newcommand{\bpf}{\begin{proof}} \newcommand{\epf}{\end{proof}}

\newcommand{\F}{\mathbb F}

\newcommand{\sA} {{\mathcal A}}
\newcommand{\sB} {{\mathcal B}}

\newcommand{\sD} {{\mathcal D}}

\newcommand{\sS} {{\mathcal S}}

\renewcommand{\phi}{\varphi}
\renewcommand{\theta}{\vartheta}

\newcommand{\coz}{{\rm coz}}
\newcommand{\wt}{{\rm wt}}

\newcommand{\mkp}{\medskip}
\newcommand{\bkp}{\bigskip}
\def\defi{\buildrel\rm def \over=}

\begin{document}

\title[Weight-preserving isomorphisms]{Weight-preserving isomorphisms between spaces of continuous functions: The scalar case}

\author[M. Ferrer]{Marita Ferrer}
\address{Universitat Jaume I, Instituto de Matem\'aticas de Castell\'on,
Campus de Riu Sec, 12071 Castell\'{o}n, Spain.}
\email{mferrer@mat.uji.es}

\author[M. Gary]{Margarita Gary}
\address{Departmento de Matem\'aticas, Universidad Aut\'onoma Metropolitana, Iztapalapa, M\'exico DF, M\'exico}
\email{gary@uji.es}

\author[S. Hern\'andez]{Salvador Hern\'andez}
\address{Universitat Jaume I, INIT and Departamento de Matem\'{a}ticas,
Campus de Riu Sec, 12071 Castell\'{o}n, Spain.}
\email{hernande@mat.uji.es}

\thanks{ The first and third listed
authors acknowledge partial support by the Generalitat Valenciana,
grant code: PROMETEO/2014/062; and by Universitat Jaume I, grant P1·1B2012-05}

\begin{abstract}
Let $\mathbb F$ be a finite field and let  $\mathcal A$  and $\mathcal B$  be vector spaces  of
$\mathbb F$-valued continuous functions  defined on locally compact spaces $X$ and $Y$, respectively.
We look at the representation of  linear bijections $H:\mathcal A\longrightarrow \mathcal B$ by
continuous functions $h:Y\longrightarrow X$ as weighted composition operators.
In order to do it, we extend the notion of Hamming metric to infinite spaces.
Our main result establishes that under some mild conditions,
every Hamming isometry can be represented as a weighted composition operator.
Connections to coding theory are also highlighted.
\end{abstract}



\date{27 January 2015}

\maketitle \setlength{\baselineskip}{24pt}

\section {Introduction}
In this paper, we are concerned with the representation of linear isomorphisms defined 
on spaces of continuous functions taking values in a vector space
$\F^n$ over a finite field $\F$. The starting point, and our main motivation,
stems from two very celebrated, and apparently disconnected,
results, whose formulation is strikingly similar, namely: MacWilliams Equivalence Theorem and 
Banach-Stone Theorem. The former one 
completely describes the isometries between block codes
(see \cite{McW:i,McW:ii}). For the reader's sake, we recall its
main features here.

Let $\F$ be a finite field. Two linear codes $C_1$ and $C_2$ over $\F$ of length $n$
are \emph{equivalent} if there is a monomial transformation $H$ of $\F^n$ such that
$T(C_1)=C_2$. Here, a monomial transformation is a linear isomorphism $H$
of the form
$$H(a_1, . . . , a_n) = (a_{\sigma(1)}w_1, . . . , a_{\sigma(n)}w_n),\ (a_1, . . . , a_n)\in \F^n,$$

\noindent where $\sigma$ is a permutation of $\{1, 2, . . . , n\}$ and $(w_1, . . . , w_n)\in(\F\setminus \{0\})^n$.
\mkp

The Hamming weight $\wt(x)$ of a vector $x\in\F^n$ is defined as the number of coordinates that are different from zero.
The following classical result establishes the relation between Hamming isometries and equivalent codes.

\bthm[MacWilliams]
Two linear codes $C_1$, $C_2$ of dimension $k$ in $\F^n$ are equivalent if and only if
there exists an abstract $\F$-linear isomorphism
$f: C_1\longrightarrow C_2$ which preserves weights, $\wt(f(x)) = \wt(x)$, for all $x\in C_1$.
\ethm
\medskip

Hence, two block codes are \emph{isometric} if
and only if they are monomially equivalent. More precisely, weight-preserving isomorphisms
between codes are given by a permutation and rescaling of the coordinates.

This fundamental result has been extended in different directions by many workers
(cf. \cite{BogartEtAlters:1978,Din_LP:2004,WardWood:1996,wood:2001}).
In particular, Heide Gluesing-Luerssen has established a variant of MacWilliams theorem for $1$-dimensional
convolutional codes and the isometries defined between them that respect the module structure of the codes
(see \cite{GluLue:2009}). It remains open the representation of general $\F$-isometries defined between convolutional
codes (cf. \cite{GluLue:2009} and \cite[Ch. 8]{piret:1988}).

The second result we are concerned in this paper, the Banach-Stone Theorem, establishes that every linear isometry
defined between the spaces of continuous functions of two compact spaces is a weighted composition operator.
It has now become a classical result that has been extended in many ways
(cf. \cite{Banach,Stone}).

\bthm[Banach-Stone Theorem]
Let $X$ and $Y$ be compact spaces and let
$H:~C(X)\longrightarrow C(Y)$ be a linear isometry.
Then $X$ and $Y$ are homeomorphic
and the isometry $H$ has the following form: there is a homeomorphism
$h: Y\longrightarrow X$, and a scalar-valued continuous function $w$ on $C(Y)$
such that  $$Hf(y)=w(y)f(h(y)),\ \forall f\in C(X),\ \forall y\in Y.$$
\ethm
\mkp

The analogy between MacWilliams and Banach-Stone theorems is blatant and
our motivation has been to explore the application of functional analysis
methods in order to extend MacWilliams Equivalence Theorem to a more general setting.
We are also concerned with the application of these
techniques to describe $\F$-isomorphisms defined between (possibly multi-dimensional) convolutional codes.

For the sake of simplicity, even though many of our results hold for spaces of group-valued continuous functions,
we shall only deal with vector-valued continuous functions on a finite field along this paper
(see \cite{FerGarHer:JFS}).
\mkp

\label{s:intro}
Let $X$ be a $0$-dimensional locally compact space,  equipped with a
Borel regular, strictly positive, measure $\mu$, and let $C_{00}(X,\F^n)$ designate the space of $\F$-valued,
compactly supported, continuous functions defined on $X$.
For any $f\in C_{00}(X,\F^n)$ and $x\in X$, we define $$\wt(f(x))\defi |\{j: \pi_j(f(x))\not= 0\}| $$ and
$$\wt(f)\defi \int_X \wt(f(x)) d\mu(x).$$
\mkp

(Notice that this integral is finite because $\wt(f(x))$ is continuous and has compact support).

The map
$$d(f,g)\defi \wt(f-g)$$
defines a metric on the vector space $C_{00}(X,\F^n)$ that is compatible with its additive group structure.
Since this metric extends the well known distance introduced by Hamming in coding theory, we  call
it \emph{Hamming metric}.\mkp

\bdfn
Let $\mathcal A$ and $\mathcal B$ be vector subspaces of $C_{00}(X,\F^n)$ and
$C_{00}(Y,\F^n)$, respectively, and let $H:\sA \longrightarrow \sB$ be a linear map.\medskip

$H$ is called \emph{Hamming isometry} if it is a linear isomorphism
and $\wt(f)=\wt(Hf)$ for each $f\in \sA$. \mkp

It is said that $H$ is a \emph{weighted composition operator} when
there exist continuous functions $h:Y\longrightarrow X$ and $w:Y\longrightarrow \F$ such that
$Hf(y)=w(y)f(h(y))$ for all $y\in Y$ and $f\in \sA$.
\edfn
\mkp

The main question we address in this research is as follows:
 \mkp

\bqtn
Is every Hamming isometry $H:\sA \longrightarrow \sB$ representable as a weighted composition
operator?
\eqtn
\mkp

In this paper, we deal with scalar-valued functions. The case of vector-valued functions will be
considered in a subsequent paper. We now introduce some pertinent notions and terminology.
\mkp

All spaces are assumed to be $0$-dimensional and Hausdorff and
throughout this paper the symbol $\F$ denotes a discrete field. 
If $X$ is a locally compact space, then $X^*$ denotes the \emph{Alexandroff compactification} of $X$,
that is, $X^*=X\cup\{\infty\}$, being $\infty$ an ideal point.

For $f\in C(X,\F^n)$, set $$\coz(f)\defi\{x\in X : f(x)\neq 0\}.$$ 
Since $\F^n$ is discrete $\coz(f)$ and $Z(f)=X\setminus \coz(f)$ are open
and closed (clopen) subsets of $X$.

Let $\sA$ be a linear subspace of $C_{00}(X,\F^n)$. For $x\in X$, let $\delta_x\colon\sA\to \F^n$ be the canonical
\emph{evaluation map} $$\delta_x(f)\defi f(x)\ \ \forall f\in\sA.$$ and
$$I_x\defi \{f\in\sA : f(x)=0\}.$$
Set  $$S\defi \{x\in X: I_x\neq \sA\}=\bigcup_{f\in\sA}\coz(f).$$

Therefore $S$ is an open subset of $X$ and, as a consequence, is also a locally compact space when
it is equipped with the topology inherited from $X$. Hence we assume WLOG that $S=X$ throughout this paper.
Thus, for each linear subspace of continuous functions considered along this paper, it is assumed:
\begin{equation}\label{(1)}
\hbox{for every}\ x\in X\ \hbox{there exists}\ f\in\sA\ \hbox{such that}\ f(x)\neq 0.
\end{equation}
\smallskip


\mkp

Define $Z(\sA)\defi \{Z(f) : f\in \sA \}$, $\coz(\sA)\defi \{\coz(f) : f\in \sA \}$, and
let $\sD$ denote the smallest ring (with respect to finite unions and intersections) of subsets containing $\coz(\sA)$.
\mkp

In coding theory, it is said that a convolutional code is \emph{controllable} when
any code sequence can be reached from the zero sequence in a finite interval
(see \cite{Fer_Her_Sha:Products,forney_trott:04,rosenthal,willems:86}).
The gist of controllability can be conveyed in a natural way to subspaces of continuous functions defined on a topological space.
In an informal way, let us say that a vector subspace of continuous functions is controllable when
any continuous functions can be reached from the zero function modulo a relatively compact open subset. It turns out
that this notion is an essential ingredient in the approach we have taken in this paper.
\mkp

\bdfn
We say that $\sA$ is \emph{controllable}\label{NC}  if
for every $f\in\mathcal{A}$ and $D_1,D_2\in\sD$ with $D_1\cap D_2=\emptyset$,
there exist $f'\in \mathcal{A}$ and  $U\in \sD$ such that
$$D_1\subseteq U\subseteq X\setminus D_2,\ f_{|D_1}=f'_{|D_1},\ \hbox{and}\  f'_{|(Z(f)\cup (X\setminus U))}=0.$$
\mkp

We say that $\sA$ \emph{separates} the points $x_1,x_2\in X$, if there is $f\in\sA$ such that
$x_1\in \coz(f)$ and $x_2\in Z(f)$ or vice versa.
\edfn
\mkp

We now formulate the main result in this paper.
\mkp

\bthm\label{mainthm}
Let $\sA$ and $\sB$ two vector spaces of\, $\F$-valued, compactly supported, continuous functions defined
on locally compact spaces $X$ and $Y$, respectively. If $\sA$ is controllable,
then every Hamming isometry $H:\sA\longrightarrow \sB$ is a weighted composition
operator.
\ethm
\mkp

\section{Basic notions and facts}

In this section, we introduce some topological notions that will be needed in the rest of the paper.
Some basic properties connecting  them are also established.

\bdfn
Two points $x_1$ and $x_2$ in $X$ are \emph{related},
written $x_1\sim x_2$, if for every $f\in \sA$ with $f(x_1)\cdot f(x_2)=0$, it follows that $f(x_1)=f(x_2)=0$. 
Let $\widetilde{X}$ be the set of equivalence classes $X/\sim$ \ equipped with the quotient topology inherited from $X$.
Every element $\widetilde{x}\in \widetilde{X}$ is associated to the coset subset $[x]\subseteq X$
consisting of all elements related to $x$. For simplicity's sake, we shall use the same symbol $[x]$ to denote
either the coset $[x]$ or the element $\widetilde x\in \widetilde X$.
Remark that $I_{x_1}=I_{x_2}$ for every $x_1$ and $x_2$ belonging to the same coset.
\edfn

\bprp\label{pesoclase}
Let $[x]$ be an equivalence class in $X$ and let $x_1,x_2\in [x]$. Then there is a unique element $\lambda(x_1,x_2)\in\F\setminus\{0\}$
such that $f(x_1)=\lambda(x_1,x_2)f(x_2)$ for all $f\in \sA$.
\eprp
\bpf
We know that $\sA\setminus I_{x}\neq \emptyset$ by (1). On the other hand, if  $f\in\sA\setminus I_{x}$, it follows that $[x]\subseteq \coz(f)$.
Pick out $x_1,x_2\in[x]$. Since $f(x_1)=f(x_1)f(x_2)^{-1}f(x_2)$, we define
\[
\lambda_f(x_1,x_2)=f(x_1)f(x_2)^{-1},
\]
which yields $f(x_1)=\lambda_f(x_1,x_2)f(x_2)$. It will suffice to verify that $\lambda_f(x_1,x_2)$
does not depend on the selected $f$  in $\sA\setminus I_{x}$.
Indeed, let $g\in\sA\setminus I_{x}$. Then $g(x_1)=\lambda_g(x_1,x_2)g(x_2)$.
The map $h\defi f(x_2)^{-1}f-g(x_2)^{-1}g\in \sA$ and $h(x_2)=0$.
Therefore $[x]\subseteq Z(h)$ and
\begin{eqnarray*}
 & 0&= \ h(x_1)\\
&  & =\ f(x_2)^{-1}f(x_1)-g(x_2)^{-1}g(x_1)\\
&  &= \ f(x_2)^{-1}\lambda_f(x_1,x_2)f(x_2) -g(x_2)^{-1}\lambda_g(x_1,x_2)g(x_2)\\
& &=\ \lambda_f(x_1,x_2)-\lambda_g(x_1,x_2).
\end{eqnarray*}

As a consequence
\[
\lambda_f(x_1,x_2)=\lambda_g(x_1,x_2)=\lambda(x_1,x_2)\in\F\setminus\{0\}.
\]
\epf
\mkp

It is readily seen that the map $\lambda(\,,\,)$ has the following properties:
\mkp

\begin{itemize}
\item$\lambda(x_2,x_1)=\lambda(x_1,x_2)^{-1},$
\item $\lambda(x_1,x_2)=\lambda(x_1,x)\lambda(x,x_2).$
\end{itemize}
\mkp

\blem\label{distintocero}
If $x_1,x_2\in X$ and $x_1 \not\sim x_2$, then there is $f_{x_1x_2}$ such that $x_1\in \coz(f_{x_1x_2})$ and $x_2\in Z(f_{x_1x_2})$.
\elem
\bpf
Since $x_1\not\sim x_2$ there is $f\in\sA$ such that $f(x_1)f(x_2)=0$ and $f(x_1)\neq 0$ or $f(x_2)\neq 0$.
If $f(x_1)\neq 0$ and $f(x_2)=0$, then $f_{x_1x_2}=f$ and we are done. Otherwise, by (1), there is $g\in\sA$ such that $g(x_1)\neq 0$.
Set $h\defi g(x_2)f-f(x_2)g\in\sA$. Then $h(x_2)=0$ and $h(x_1)=-f(x_2)g(x_1)\neq 0$. In this case $f_{x_1x_2}=h$.
\epf

\bdfn
$A\subseteq X$ is called \emph{saturated} if and only if $x\in A$ implies $[x]\subseteq A$.
\edfn
\mkp

The proof of the next result is easy. We include it for the sake of completeness.

\mkp


\bprp\label{compacto}
For every $f\in \sA$ and $x\in X$, we have:
\begin{enumerate}
\item[(a)] $\coz(f)$ and $Z(f)$ are saturated subsets of $X$.
\item[(b)] $[x]$ is a saturated compact subset of $X$.
\end{enumerate}
\eprp
\bpf
The proof of (a) is clear. (b) Let $x\in X$. We first proof that $[x]$ is closed in $X$. Let $x'\in X\setminus [x]$.
By Lemma \ref{distintocero} there is $f\in\sA$ such that $x'\in \coz(f)$ and $x\in Z(f)$.
Applying (a), it follows that $[x']\subseteq \coz(f)$ and $[x]\subseteq Z(f)$.
 Then $x'\in \coz(f)\subseteq X\setminus [x]$ and $\coz(f)$ is open in $X$.

On the other hand, by (1), there is $g\in\sA$ such that $[x]\subseteq \coz(g)$.
Since $\coz(g)$ is compact and $[x]$ is closed in $X$, we have that $[x]$ is compact.
\epf
\mkp

 Let $\pi\colon X\to \widetilde{X}$ denote the canonical quotient map associated to the equivalence relation $\sim$
 and equip  $\widetilde{X}$ with the canonical quotient topology.
Using Proposition \ref{compacto}, it is easily seen that the subsets $\pi(\coz(f))$ and $\pi(Z(f))$ are clopen in $\widetilde{X}$
for every $f\in\sA$ and, with a little more effort, it is proved that $\widetilde X$ is a Hausdorff, locally compact space.
We leave the verification of this fact to the interested reader.
\mkp


A standard compactness argument is used in the proof of the following lemma. We include it here for the sake of completeness.

\blem\label{lema0}
Let $K_1$ and $K_2$ be compact subsets of $X$ such that $x_1\not\sim x_2$  for every $x_1\in K_1$ and $x_2\in K_2$.
Then there are $D_1,D_2\in\sD$
such that $K_1\subseteq D_1$, $K_2\subseteq D_2$ and $D_1\cap D_2=\emptyset$.
\elem
\bpf
Let $x_1\in K_1$ and $x\in K_2$, which implies $x_1\not\sim x$.
By Lemma \ref{distintocero}, there is $f_x\in\sA$ such that $[x_1]\subseteq \coz(f_x)$ and $[x]\subseteq Z(f_x)$.
We have $K_2\subseteq\bigcup\limits_{[x]\in\pi(K_2)} Z(f_x)$ and $[x_1]\subseteq\bigcap\limits_{[x]\in\pi(K_2)}\coz(f_x)$. Since $K_2$ is compact and $Z(f_x)$ is open,
we have $K_2\subseteq\bigcup\limits_{i=1}^n Z(f_{x^{(i)}})$ and
$[x_1]\subseteq \bigcap\limits_{i=1}^n\coz(f_{x^{(i)}})=X\setminus\bigcup\limits_{i=1}^nZ(f_{x^{(i)}})\subseteq X\setminus K_2$.

Define $C_{x_1}=\bigcap\limits_{i=1}^n \coz(f_{x^{(i)}})$, which is a clopen subset of $X$. Remark that $[x_1]\subseteq C_{x_1}$ and
$C_{x_1}\cap K_2=\emptyset$. Consequently $K_1\subseteq\bigcup\limits_{[x]\in\pi(K_1)}C_x$ and $C_x\cap K_2=\emptyset$ for every $[x]\in\pi(K_1)$.
Since $K_1$ is compact, we have $K_1\subseteq \bigcup\limits_{j=1}^m C_{x_{(j)}}$.

Define $D_1=\bigcup\limits_{j=1}^mC_{x_{(j)}}\in\sD$ and observe that $K_1\subseteq D_1$ and $D_1\cap K_2=\emptyset$.
Since $D_1$ is a saturated compact subset of $X$, we repeat again the same procedure in order to obtain
$D_2\in\sD$ such that $K_2\subseteq D_2$ and $D_1\cap D_2=\emptyset$.
\epf
\mkp

We notice that the lemma above applies to any two disjoint saturated compact subsets of $X$. On the other hand,
the following remark is easily seen.

\brem\label{abiertos} Every $D\in\sD$ is a saturated compact subset of $X$ and $\pi(D)$ is clopen in $\widetilde{X}$.
Furthermore, the collection
$\{\pi(D) : D\in \sD \}$ is an open base for $\widetilde X$.
\erem

\section{Separating maps and support subsets}
\bdfn
A map $H:\sA\longrightarrow \sB$ is said to be \emph{separating} (or \emph{disjointness preserving})
when $\coz(f)\cap\coz(g)=\emptyset$ implies $\coz(Hf)\cap\coz(Hg)=\emptyset$, $f,g\in \sA$. 
\mkp

A linear functional $\phi:\sA\longrightarrow \F$ is called \emph{separating} when $\coz(f)\cap\coz(g)=\emptyset$
implies $\phi(f)\cdot \phi(g)=0$.
\edfn
\mkp

\blem
Let $f$ and $g$ be two elements in $\sA$. Then $\coz(f)\cap\coz(g)=\emptyset$ if and only if
$\wt(f+g)=\wt(f)+\wt(g)$.
\elem
\bpf
It follows from the inequality
$$\wt(f+g)\leq \wt(f)+\wt(g)-\wt(f\cdot g)$$
that is readily verified.
\epf
\mkp

\bcor\label{hammingcor}
Every Hamming isometry is a separating linear isomorphism.
\ecor


Separating isomorphisms have been studied by many workers and have found application to a variety of fields
(cf. \cite{Abram_Kito:2000,Alaminos_et_Alters:2009,Ara_Bec_Nar,Ara_Jar:2003,Burgos_et_Alters:2013,Che_Ke_Lee_Won:2003,Chebotar_et_Alters:2004,
font-her,font-her2,gau-jeang-wong,her-bec-nar,jarosz,Lau_Wong:2013}).
After Corollary \ref{hammingcor}, it is clear that, in order to prove Theorem \ref{mainthm},
it suffices to deal with the broader case of separating isomorphisms and
so we do in the rest of the paper.
\mkp

The following definition makes sense for every subset of $ X $ but we have restricted it to saturated subsets,
because it will only be applied to these subsets in this paper. \mkp

\bdfn
Let $\phi :\sA\longrightarrow \F$ be a map. A saturated closed subset $K$ of $X$ is said to be a \emph{support}
for $\phi$ if given $f\in\sA$ with $K\subseteq Z(f)$, it holds that $\phi(f)=0$.
\edfn
\mkp

Support subsets enjoy several nice properties.

\bprp\label{propiedadessoporte}
Let $\phi :\sA\longrightarrow \F$ be a non null, separating, linear functional.
Then the following assertions hold:
\begin{enumerate}
\item[(a)] $X$ is a support for $\phi$.
\item[(b)] If $K$ is a support for $\phi$ then $K\neq\emptyset$.
\item[(c)] Let $K$ be a support for $\phi$ and $f,g\in\sA$ such that $f_{|K}=g_{|K}$.
Then $\phi(f)=\phi(g)$.
\item[(d)] If $\sA$ is controllable and $K_1$ and $K_2$ are both supports for $\phi$, then $K_1\cap K_2\neq\emptyset$.
\end{enumerate}
\eprp
\begin{proof}
(a) This is clear. 

\noindent (b) Let $K$ be a support for $\phi$ and suppose $K=\emptyset$. Then $K=\emptyset\subseteq Z(f)$ for all $f\in\sA$. Consequently
$\phi(f)=0$ for all $f\in\sA$, which is a contradiction since $\phi$ is non null.

\noindent (c) Let $K$ be a support for $\phi$. If $f,g\in\sA$ and $f_{|K}=g_{|K}$ then $f-g\in\sA$ and $K\subseteq Z(f-g)$.
So $0=\phi(f-g)=\phi(f)-\phi(g)$.

\noindent (d) Let $K_1$ and $K_2$ be supports for $\phi$ and suppose that $K_1\cap K_2=\emptyset$.
Since $\phi$ is non null, there is $f\in\sA$ such that $\phi(f)\neq 0$. Remark that the set $C_1=\coz(f)\cap K_1\neq\emptyset$ because,
otherwise, $K_1\subseteq Z(f)$ and then $\phi(f)=0$, which is not true. Since $\coz(f)$ is a saturated compact subset of $X$ and
$K_1$ is also saturated and closed, it follows that $C_1$ is a saturated compact subset of $X$. In like manner
$C_2=\coz(f)\cap K_2$ is non empty, saturated and compact. Furthermore $C_1\cap C_2=\emptyset$ and by Lemma \ref{lema0} there exist
$D_1,D_2\in\sD$ such that $C_1\subseteq D_1$, $C_2\subseteq D_2$ and $D_1\cap D_2=\emptyset$.
Applying that $\sA$ is controllable to $D_1$, $D_2$ and $f$, we obtain  $U\in \sD$ and $f'\in\mathcal A$ such that
$C_1\subseteq D_1\subseteq U\subseteq X\setminus D_2\subseteq X\setminus C_2$ and
$f_{|D_1}=f'_{|D_1}\ \hbox{and}\ \  f'_{|(Z(f)\cup (X\setminus U))}=0.$

Remark that $\coz(f)=C_1\cup C_2\cup(\coz(f)\setminus(C_1\cup C_2))$.  Evaluating $f'$ yields:
\begin{enumerate}
\item [] If $x\in C_1$ then $f'(x)=f(x)$.
\item [] If $x\in K_1\setminus C_1$ then $f'(x)=0=f(x)$.
\item [] If $x\in K_2$ then $f'(x)=0$.
\end{enumerate}

As a consequence $f'_{|K_1}=f_{|K_1}$ and $f'_{|K_2}=0$.
Applying Proposition \ref{propiedadessoporte}, we deduce that $\phi(f')=\phi(f)\neq 0$ and $\phi(f')=0$, which
is a contradiction. This completes the proof.

\end{proof}
\mkp

Next it is proved that, when $\sA$ is controllable, every non null, separating, linear functional $\phi:\sA\longrightarrow \F$ has
a minimum support set. For that purpose,
we define 
$$\sS=\{A\subseteq X : A \text{ is support for }\phi \}.$$

There is a canonical partial order that can be defined on $\sS$:
$A\leq B$, $A,B\in \sS$, if and only if $B\subseteq A$. A standard compactness argument shows that
$(\sS,\leq)$ is an inductive set and, by Zorn's lemma, $\sS$ has a $\subseteq$-minimal element $K$.
\mkp

\bprp\label{1elemento}
Let $\phi :\sA\longrightarrow \F$ be a non null, separating, linear functional. If $\sA$ is controllable,
then there exists $x\in X$ such that
$K=[x]$  is a support for $\phi$.
\eprp
\bpf
By Proposition \ref{propiedadessoporte} $K\neq\emptyset$. Suppose now that there are two different cosets $[x_1],[x_2]$
that are contained in $K$. Since $X$ is Hausdorff and $K$ is saturated, using  Lemma \ref{lema0},
we can select two disjoint
saturated open sets  $V_1,V_2\subseteq X$ such that $[x_1]\subseteq V_1$ and $[x_2]\subseteq V_2$. Since $K$ is minimal, the subset
$K\setminus V_i$ is a saturated closed subset of $X$ that is not a support for $\phi$. Hence, there is $f_i\in\sA$
such  that $K\setminus V_i\subseteq Z(f_i)$ and $\phi(f_i)\neq 0$, $1\leq i\leq 2$. As $\phi$ is a separating functional,
the subset
$A=\coz(f_1)\cap \coz(f_2)$ is a nonempty saturated compact subset of $X$. We claim that
$K\cap A=\emptyset$. Indeed, otherwise, pick out an element $a\in K\cap A$.
Then $[a]\subseteq  K\cap A$. If $[a]\subseteq V_1$ then $[a]\subseteq K\setminus V_2$ and
$[a]\subseteq Z(f_2)$, which is a contradiction. On the other hand, if $[a]\nsubseteq V_1$ then $[a]\subseteq K\setminus V_1$
and $[a]\subseteq Z(f_1)$, which is a contradiction again. Therefore, we have proved that $K\cap A=\emptyset$.

Take now $B=K\cap(\coz(f_1)\cup \coz(f_2))$. If $B=\emptyset$ then $K\cap \coz(f_i)=\emptyset$ and $K\subseteq Z(f_i)$,
which implies $\phi(f_i)=0$, $1\leq i\leq 2$, and we obtain a contradiction. Therefore, we have $B\neq\emptyset$.
Thus $B$ is a saturated compact
subset of $X$ satisfying that $A\cap B=\emptyset$. Applying Lemma \ref{lema0}, we can select two disjoint subsets
$D_A,D_B\in\sD$ such that $A\subseteq D_A$ and $B\subseteq D_B$. Applying that $\sA$ is controllable to $D_A$, $D_B$ and $f_1$,
we can take $U\in \sD$ and $f'\in \sA$ such that
$B\subseteq D_B\subseteq U\subseteq X\setminus D_A\subseteq X\setminus A$, which implies $U\cap A=\emptyset$,
$f_{1|D_B}=f'_{|D_B}\ \hbox{and}\ \  f'_{|(Z(f_1)\cup (X\setminus U))}=0.$

Let us see that $f'_{|K}= f_{1|K}$. Indeed, if $x\in K\setminus \coz(f_1)$ then $f'(x)=0=f_1(x)$ and
if $x\in K\cap \coz(f_1)\subseteq D_B$ then $f'(x)=f_1(x)\neq 0$. By Proposition \ref{propiedadessoporte}
$\phi(f')=\phi(f_1)\neq 0$. Since $\phi$ is separating,
$\emptyset\neq \coz(f')\cap \coz(f_2)\subseteq \coz(f_1)\cap \coz(f_2)=A$. But this is a contradiction because
$A\subseteq Z(f')$. By Proposition \ref{compacto}, it follows that $K$ may only contain an equivalence class $[x]=K$,
for some point $x$ in $X$. This completes the proof.
\epf
\bkp

\section{Proof of main result}

We have remarked after Corollary \ref{hammingcor} that, in order to prove the main result formulated at the Introduction,
it suffices to deal with separating linear isomorphisms. Therefore, assume that $H:\sA\longrightarrow \sB$ is a
linear separating map defined between linear subspaces $\sA$ and $\sB$ of $C_{00}(X,\F)$ and $C_{00}(Y,\F)$, respectively.
Observe that for every $y\in Y$, the composition $\delta_y\circ H$ is a separating linear functional of
$\sA$ into $\F$. Conveying to $Y$ and $\sB$ the equivalence relation we have defined above on $X$ and $\sA$, and
applying to $\delta_{y}\circ H$ the last two results in the previous section, we obtain:
\mkp

\bprp
Let $H:\sA\longrightarrow \sB$ be a linear separating map. If  $K$ is a support for $\delta_y\circ H$ and $y'\in [y]$
then $K$ is a support to $\delta_{y'}\circ H$.
\eprp
\bpf
It suffices to take into account that every $Z\in Z(\sB)$ is saturated.
\epf
\mkp

Applying Proposition \ref{1elemento} to $\delta_y\circ H$, for each $y\in Y$,
we are now in position of defining the \emph{support map} $h$ that is associated to $H$.
This map is defined between the spaces $Y$ and  $\widetilde X$. Again, in order to
simplify the notation, we will use the same symbol $h(y)$ to denote both, an element of $\widetilde X$,
and the equivalence class $\pi^{-1}(h(y))$, which is a subset of $X$.
\mkp

\bprp\label{propiedadesh}
Let $H:\sA\longrightarrow \sB$ a separating linear map satisfying that for every $y\in Y$ there is
$f_y\in \sA$ such that $Hf_y(y)\neq 0$. If $\sA$ is controllable, then there is a map
$h : Y\longrightarrow \widetilde{X}$ 
satisfying the following properties:\medskip

\item [(a)]  For every $f\in\sA$ with $f_{|h(y)}=0$, it follows that $Hf(y)=0$.
\item [(b)] $h(y')=h(y)$ for all $y'\sim y$.
\item [(c)] If $A\subsetneq\widetilde{X}$ is open, $f\in\sA$ and $\pi^{-1}(A)\subseteq Z(f)$ then $h^{-1}(A)\subseteq Z(Hf)$.
\item [(d)] $h(\coz(Hf))\subseteq \pi(\coz(f))$ for every $f\in \sA$.
\eprp
\bpf
We define $h(y)$ as the smallest support associated to $\delta_y\circ H$.

(a) This is clear.

(b) It follows from $\sS_y=\sS_{y'}$ when $y\sim y'$.

(c) Take $y\in h^{-1}(A)$. Then $\pi^{-1}(\widetilde{X}\setminus A)$ is a nonempty, saturated, and closed subset
that it is not a support for $\delta_y\circ H$. Therefore, there is $g\in\sA$ such that $\pi^{-1}(\widetilde{X}\setminus A)\subseteq Z(g)$
and $Hg(y)\neq 0$. So we have $\coz(g)\subseteq \pi^{-1}(A)$ and $\coz(f)\subseteq X\setminus\pi^{-1}(A)$.
Since $H$ is a separating map, $\coz(Hg)\cap \coz(Hf)=\emptyset$. As a consequence $Hf(y)=0$.

(d) Let $[x]\in h(\coz(Hf))$, then $[x]=h(y)$ for some $y\in \coz(Hf)$. Since $h(y)$ is support for $\delta_y\circ H$,
we have $[x]\nsubseteq Z(f)$. Since $Z(f)$ is saturated, it follows that $[x]\subseteq \coz(f)$.
\epf
\mkp

Let $Gr[h]\defi \bigcup\limits_{y\in Y} (h(y)\times \{y\})$ denote the graphic of $h$
equipped with the topology inherited as a subspace of
$X\times Y$. We have the following representation of separating linear maps. \medskip

\bprp\label{representaciondeH}
Let $H:\sA\longrightarrow \sB$ a separating linear map satisfying that for every $y\in Y$ there is
$f_y\in \sA$ such that $Hf_y(y)\neq 0$. If $\sA$ is controllable, then there is a  map
$\omega :Gr[h]\longrightarrow \F\setminus\{0\}$
satisfying the following properties: \medskip

\begin{enumerate}
\item [(a)] $Hf(y)=\omega(x,y)f(x)$ for all $(x,y)\in Gr[h]$ and all $f\in\sA$.
\item [(b)] $\omega(x',y')=\lambda(y',y)\omega(x,y)\lambda(x,x')$ for all $y'\sim y$ and $(x,y) ,(x',y')\in Gr[h]$.
\item [(c)] $\omega$ is continuous.
\end{enumerate}
\eprp
\bpf
(a) Let $(x,y)\in Gr[h]$. By hypothesis, there is $f'\in \sA$ such that $Hf'(y)\neq0$.
Then $f'(x)\neq0$ since $h(y)$ is a support set for $\delta_y\circ H$.
Set $\alpha=f'(x)\in\F\setminus\{0\}$ and $f_x=\alpha^{-1}f'\in\sA$, which implies $f_x(x)=1$.

We define $$\omega(x,y)=Hf_x(y)=\alpha^{-1}Hf'(y)\in\F\setminus\{0\}.$$

Observe that $\omega(x,y)$ does not depend on the specific map $f\in\sA$ with $f(x)=~1$ we select.
Indeed, let $g_x\in\sA$ such that $g_x(x)=1$. Take $x'\in h(y)$, then by Proposition \ref{pesoclase}
$f_x(x')=\lambda(x',x)f_x(x)=\lambda(x',x)=\lambda(x',x)g_x(x)=g_x(x')$. Thus, we have shown that
$(f_x)_{|h(y)}=(g_x)_{|h(y)}$. By Proposition \ref{propiedadessoporte}, we have $Hg_x(y)=Hf_x(y)=\omega(x,y)$.

Pick out now an arbitrary map $f\in\sA$. If $f(x)=0$ then, since $Z(f)$ is saturated, $h(y)=[x]\subseteq Z(f)$
and $Hf(y)=0$. Obviously $Hf(y)=\omega(x,y)f(x)=0$. Therefore,
suppose WLOG that $f(x)=\beta\neq0$ and set $g'_x=\beta^{-1}f\in\sA$. Then we have $g'_x(x)=1$ and,
since $\omega(x,y)$ does not depend on $g'_x$, it follows that
$Hg'_x(y)=Hf_x(y)=\omega(x,y)$. Taking into account that $H$ is a linear map, we get
$Hg'_x(y)=\beta^{-1}Hf$. Thus $\beta^{-1}Hf(y)=\omega(x,y)$, which yields
$Hf(y)=\beta\omega(x,y)=\omega(x,y)f(x)$. This completes the proof. \mkp

(b) This is clear after making some straightforward evaluations.\mkp

(c) Let $((x_d,y_d))_d$ be a net converging to $(x,y)$ in $Gr[h]$ and  
 take $f_x\in \sA$ such that $f_x(x)=1$. Since $\F$ is discrete and $f_x$ and $Hf_x$ are continuous,
 there exists $d_0$ such that $f_x(x_d)=1$ and $Hf_x(y_d)=Hf_x(y)$ for all $d\geq d_0$.
 Thus $\omega(x_d,y_d)=\omega(x_d,y_d)f_x(x_d)=Hf_x(y_d)=Hf_x(y)=\omega(x,y)f_x(x)=\omega(x,y)$ for all
 $d\geq d_0$. This implies that the net $(\omega(x_d,y_d))_d$ converges to $\omega(x,y)$.
\epf
\mkp

As a consequence of the previous result, we obtain a converse to Proposition \ref{propiedadesh}.

\bcor\label{consecuenciarepresentacion}
$Hf(y)=0$ implies $f(x)=0$ for all  $(x,y)\in Gr[h]$ .
\ecor
\mkp

Our next goal is to verify that the support map $h$ is continuous and surjective assuming the same conditions as
in Proposition \ref{propiedadesh} if $H$ is also one-to-one.
We split the proof in several lemmata for the reader's sake.
\mkp

\blem\label{continuidadh}
Assuming the same conditions as in Proposition \ref{propiedadesh}, the support map $h\colon Y\to\widetilde{X}$ is continuous.
\elem
\bpf
Let  $(y_d)_{d\in D}$ be a net in $Y$ converging to $y\in Y$. Since $\widetilde{X}$ is locally compact and Hausdorff, its
Alexandroff compactification $\widetilde{X}^*$ is also Hausdorff. By a standard compactness argument, we may assume WLOG
that $(h(y_d))_d$ converges to $t\in\widetilde{X}^*$. Reasoning by contradiction, suppose $h(y)\neq t$ and
take two disjoint open neighborhoods $V_{h(y)}$ and $V_t$ of $h(y)$ and $t$ respectively.
Take $d_1$ such that $h(y_d)\in V_t \cap\widetilde{X}$ for all $d\geq d_1$.

 Since the support sets for $\delta_{z}\circ H$ contains $h(z)$ for all $z\in Y$, it follows that
  the subset $\pi^{-1}(\widetilde{X}\setminus(V_{h(y)}\cap \widetilde{X}))$ may not be a support set for $\delta_y\circ H$. Therefore,
  there exists $f\in\sA$ such that $\pi^{-1}(\widetilde{X}\setminus(V_{h(y)}\cap \widetilde{X}))\subseteq Z(f)$ and $Hf(y)\neq0$.
  Moreover, since $H(f)$ is continuous, the net $(Hf(y_d))_{d\in D}$ converges to $Hf(y)$ and, since $\F$ is discrete,
  there is $d_2\geq d_1$ such that $Hf(y_d)\neq0$ for all $d\geq d_2$. Therefore, the subset
  $\pi^{-1}(\widetilde{X}\setminus(V_t\cap\widetilde{X}))$ may not be a support set for $\delta_{y_{d_3}}\circ H$ for some index
  $d_3\geq d_2$. As a consequence, there exists $f_3\in\sA$ such that
  $\pi^{-1}(\widetilde{X}\setminus(V_{t}\cap\widetilde{X}))\subseteq Z(f_3)$ and $Hf_3(y_{d_3})\neq0$.
  Thus, we have $y_{d_3}\in \coz(Hf_3)\cap \coz(Hf)$ and, since $H$ is a separating map, $\coz(f_3)\cap \coz(f)\neq\emptyset$.
  But $\coz(f_3)\subseteq \pi^{-1}(V_t\cap\widetilde{X})$ is disjoint from  $\coz(f)\subseteq \pi^{-1}(V_{h(y)}\cap\widetilde{X})$.
  This contradiction
  completes the proof.
\epf
\mkp

\blem\label{densidadh}
Assuming the same conditions as in Proposition \ref{propiedadesh}, if $H$ is also one-to-one,
then $h(Y)$ is dense in $\widetilde{X}$.
\elem
\bpf
Reasoning by contradiction again, suppose there is $x\in X$ such that
$[x]\notin \overline{h(Y)}^{\widetilde{X}}$. Set $A=\overline{h(Y)}^{\widetilde{X}}$, which implies $[x]\cap\pi^{-1}(A)=\emptyset$.
On the other hand, by (\ref{(1)}), there is $f\in\sA$ such that $[x]\subseteq \coz(f)$. Define $B=\pi^{-1}(A)\cap \coz(f)$, which is a saturated compact
subset because $\pi^{-1}(A)$ is closed and $\coz(f)$ is compact and saturated. Moreover, we have that $B\neq\emptyset$.
Otherwise, $\pi^{-1}(h(Y))\subseteq\pi^{-1}(A)\subseteq Z(f)$. This implies that $Hf\equiv0$ and $f\equiv0$, which is a contradiction.
Since $[x]\cap B=\emptyset$, by Lemma \ref{lema0}, there are two disjoint subsets $D_x,D_B\in\sD$ such that
$[x]\subseteq D_x$ and $B\subseteq D_B$. Then the subset $D=D_x\cap \coz(f)\in\sD$ contains $[x]$ and
$D\cap\pi^{-1}(A)=\emptyset$. We now apply that $\sA$ is controllable to $D$, $D_B$ and $f$ in order to obtain
$U\in \sD$ and $f'\in \sA$
such that $[x]\subseteq D\subseteq U\subseteq X\setminus D_B\subseteq X\setminus B$,
$f_{|D}=f'_{|D}\ \hbox{and}\ \  f'_{|(Z(f)\cup (X\setminus U))}=0.$
 Hence $\coz(f')\subseteq U\cap \coz(f)$, $U\cap B=\emptyset$ and $\coz(f')\cap\pi^{-1}(A)=\emptyset$.
 As a consequence $\pi^{-1}(h(Y))\subseteq\pi^{-1}(A)\subseteq Z(f')$ and $Hf(y)=0$ for all $y\in Y$.
 Since $H$ is a linear monomorphism we have $f\equiv0$, which is a contradiction. Therefore $\overline{h(Y)}
^{\widetilde{X}}=\widetilde{X}$, which completes the proof.
\epf
\mkp

Let $Y^*$ and $\widetilde{X}^*$ be the Alexandroff compactification of $Y$ and $\widetilde{X}$ respectively. Then
there is a canonical way of extending $h$ to a map $h^*\colon Y^*\to\widetilde{X}^*$ by $h^*|_Y=h$ and $h^*(\infty)=\infty$.
It turns out that this canonical extension is a continuous onto map. \mkp

\blem\label{extensionh}
Assuming the same conditions as in Proposition \ref{propiedadesh}, if $H$ is also one-to-one,
then $h^*$ is continuous and onto.
\elem
\bpf

Since $h^*|_Y=h$ is continuous, in order to prove the continuity of  $h^*$, it suffices to verify the continuity of $h^*$ at $\infty$.
Reasoning by contradiction, suppose that $h^*$ is not continuous at $\infty$. Then, there must be a compact subset $K_0\subseteq \widetilde{X}$
such that $\infty\in\overline{h^{-1}(K_0)}^{Y^*}$. Otherwise, we would have $\infty\notin\overline{h^{-1}(K)}^{Y^*}$ for every compact subset $K$ of $\widetilde{X}$.
Since $h^{-1}(K)$ is closed in $Y$, it follows that $h^{-1}(K)=\overline{h^{-1}(K)}^Y=\overline{h^{-1}(K)}^{Y^*}$.
However, every closed subset of $Y^*$ is either the union of $\{\infty\}$ and a closed subset of $Y$, or a compact subset of $Y$.
Hence  $h^{-1}(K)$ is compact in $Y$ for every compact subset $K$ in $\widetilde{X}$ and, as a consequence, we have
 $\infty\in Y^*\setminus h^{-1}(K)$, which is open in $Y^*$. Thus, we have proved that
 $\widetilde{X}^*\setminus K$ is an open neighborhood of $\infty=h^*(\infty)$ and
 $h^*(\infty)\in h^*(Y^*\setminus h^{-1}(K))\subseteq \widetilde{X}^*\setminus K$ for every compact subset $K$ of $\widetilde{X}$,
which would yield the continuity of $h^*$ at $\infty$.

Take a net $(y_d)_{d\in D}\subseteq h^{-1}(K_0)$ converging to $\infty$. By the compactness of $K_0$, we may assume WLOG that
$(h(y_d))_{d\in D}$ converges to $[x_0]\in K_0$. But  $\coz(Hf)$ is compact and $\infty\in Y^*\setminus \coz(Hf)$ for all $f\in\sA$.
Therefore, for every $f\in \sA$, there is an index $d(f)$ such that $y_d\in Y\setminus \coz(Hf)$ for all $d\geq d(f)$. That is
$Hf(y_d)=0$  and, by Corollary \ref{consecuenciarepresentacion}, we have  $f_{|h(y_d)}= 0$ for all $d\geq d(f)$.
Thus $(h(y_d))_{d\geq d(f))}$ is contained in $\pi(Z(f))$ and, as a consequence, we have
$[x_0]\in\overline{\pi(Z(f))}^{\widetilde{X}}=\pi(Z(f))$ for all $f\in\sA$. This implies that
$f(x_0)=0$ for all $f\in \sA$, which is a contradiction.

Now, it is easy to show that $h^*$ is an onto map. Indeed,
since $Y^*$ is compact, $h^*$ is continuous and $\widetilde{X}^*$ is Hausdorff,
we have that $h^*(Y^*)$ is a compact subset of $\widetilde{X}^*$. Therefore
$\overline{h^*(Y^*)}^{\widetilde{X}^*}=h^*(Y\cup\{\infty\})=h(Y)\cup\{\infty\}\subseteq\overline{h(Y)}^{\widetilde{X}^*}\cup\{\infty\}=\overline{h^*(Y^*)}^{\widetilde{X}^*}$
and, by Lemma \ref{densidadh}, it follows that
$h^*(Y^*)=\overline{h^*(Y^*)}^{\widetilde{X}^*}=\overline{h(Y)}^{\widetilde{X}}\cup\{\infty\}=\widetilde{X}\cup\{\infty\}=\widetilde{X}^*$.
\epf
\mkp

From Proposition \ref{extensionh}, it follows a main partial result.

\bcor\label{teoremah}
Assuming the same conditions as in Proposition \ref{propiedadesh}, if $H$ is also one-to-one,
then $h\colon Y\to\widetilde{X}$ is continuous and onto.
\ecor
\mkp

Set $\widetilde{h}\colon\widetilde Y\to \widetilde X$ by $\widetilde h([y])=h(y)$ for all $[y]\in \widetilde Y$,
which is clearly well defined. A straightforward consequence of Corollary \ref{teoremah} is: \mkp

\bprp
Assuming the same conditions as in Proposition \ref{propiedadesh}, if $H$ is also a bijection,
then $\widetilde{h}$ is a homeomorphism of\ $\widetilde Y$ onto $\widetilde X$.
\eprp
 \bpf
 The continuity of $\widetilde{h}$ follows from the continuity of $h$ and $\pi$.

 Take $[y_1]\neq [y_2]$ in $Y$. By Lemma \ref{distintocero}, there is $f\in\sA$ such that
 $[y_1]\subseteq Z(Hf)$ and $[y_2]\subseteq \coz(Hf)$.
  Applying Corollary \ref{consecuenciarepresentacion} and Proposition \ref{propiedadesh}, we obtain
  $h(y_1)\subseteq Z(f)$ and $h(y_2)\subseteq \coz(f)$, which implies
 $\widetilde h([y_1])\neq \widetilde h([y_2])$. Thus $\widetilde h$ is $1$-to-$1$.
 On the other hand, the map $\widetilde h$ is onto because so is $h$.

 Now, we can proceed as in Lemma \ref{extensionh}, in order to extend $\widetilde h$ to
 a continuous map $\widetilde h^*\colon \widetilde Y^*\to\widetilde X^*$. Clearly the map
  $\widetilde h^*$ is a continuous bijection and, therefore a homeomorphism between
  compact spaces. This automatically implies that $\widetilde h$ is a homeomorphism.
  \epf
\mkp

We can now establish the representation of separating isomorphisms as weighted composition operator,
which implies Theorem \ref{mainthm}.

\bthm\label{thm:separating}
Let $H:\sA\longrightarrow \sB$ a separating, linear, onto, map. 
If $\sA$ is controllable, then there are continuous maps
$h : {Y}\longrightarrow \widetilde{X}$ and $\omega : Gr[h]\longrightarrow \F$  satisfying the following properties:
\medskip

\begin{enumerate}
\item [(a)] For each $y\in Y$, $x\in {h}(y)$, and every $f\in \sA$
it holds $$Hf(y)=\omega(x,y)f(x).$$
\item [(b)] $H$ is continuous with respect to the pointwise convergence topology.
\item [(c)] $H$ is continuous with respect to the compact open topology.

\end{enumerate}
\ethm
\bpf
Since both $\sA$ and $\sB$ satisfy the initial assumption (\ref{(1)}), it follows
that item (a) is a direct consequence from Proposition \ref{representaciondeH}.
On the other hand, it is readily seen that (a) implies (b). Thus only (c) needs verification.

(c)\ Let $(f_d)_d\subseteq\sA$ be a net uniformly converging to $0$  in the compact open topology.
If $K$ is a compact subset of $Y$, then $h(K)$ is a compact subset of $\widetilde X$ by the continuity of $h$.
Furthermore, by Remark \ref{abiertos}, the subset $\pi^{-1}(h(K))$ is compact in $X$ . Indeed, for every
$[x]\in  h(K)$, there is $f_x\in\sA$ such that $[x]\in \pi(\coz(f_x))$. Hence
$h(K)\subseteq \bigcup\limits_{[x]\in h(K)} \pi(\coz(f_x))$. By compactness, there is a finite subcover, say
$h(K)\subseteq \bigcup\limits_{1\leq i\leq n} \pi(\coz(f_i))$. Thus
$\pi^{-1}(h(K))\subseteq \bigcup\limits_{1\leq i\leq n} \coz(f_i)$, which yields the compactness of $\pi^{-1}(h(K))$.

Since $(f_d)_d$ converges to $0$ uniformly on  $\pi^{-1}(h(K))$, it follows that $(f_d)_d$ is eventually equal to $0$
on $\pi^{-1}(h(K))$. Applying (1), it follows that $(Hf_d)_d$ is eventually $0$ on $K$. This completes the proof.
\epf
\mkp

We are now in position of establishing the main result formulated at the Introduction.
\mkp

\begin{proof}[Proof of Theorem \ref{mainthm}]
Since $H$ is a Hamming isometry of $\sA$ onto $\sB$, it is separating by Corollary \ref{hammingcor}.
Thus $H$ must be a weighted composition operator by Theorem \ref{thm:separating}.
\end{proof}


\begin{thebibliography}{99}

\bibitem{Abram_Kito:2000}
{Abramovich, Y. A. and Kitover, A. K.},
{\it Inverses of disjointness preserving operators},
{Mem. Amer. Math. Soc.},
{\bf 143}, {(2000)}, {no. 679}.

\bibitem{Alaminos_et_Alters:2009}
{Alaminos, J. and Bre{\v{s}}ar, M. and Extremera, J. and Villena, A. R.},
{\it Maps preserving zero products}, {Studia Math.},
{\bf 193}, no. 2, {(2009)}, {131--159}.

\bibitem{Ara_Bec_Nar}
Araujo, J., Beckenstein, E., Narici, L.,
{\it On biseparating maps between realcompact spaces},
Journal of Mathematical Analysis and Applications,
{\bf }, no.

\bibitem{Ara_Jar:2003}
{Araujo, J. and Jarosz, K.:},
{\it Biseparating maps between operator algebras},
  {Journal of Mathematical Analysis and Applications},
   {\bf 282}, (1),
      {(2003)},
    {48--55}.

\bibitem{Banach}
Banach, S.: {\it Th\'eorie des op\'erations lin\'eaires}, Chelsea, 1955.

\bibitem{BogartEtAlters:1978}
{Bogart, K., Goldberg, D., and Gordon, J.},
     {\it An elementary proof of the {M}ac{W}illiams theorem on
              equivalence of codes},
   {Information and Computation},
   {\bf 37}, no. 1, {(1978)},
    {19--22}.

\bibitem{Burgos_et_Alters:2013}
{Burgos, Mar{\'{\i}}a and S{\'a}nchez-Ortega, Juana},
{\it On mappings preserving zero products},
{Linear Multilinear Algebra}, {\bf 61}, no. 3, {(2013)}, {323--335}.

\bibitem{Che_Ke_Lee_Won:2003}
{Chebotar, M. A., Ke, W.-F., Lee, P.-H., and Wong, N.-C.},
{\it Mappings preserving zero products},
{Studia Math.}, {\bf 155}, no. 1, {(2003)}, {77--94}.

\bibitem{Chebotar_et_Alters:2004}
{Chebotar, M. A. and Ke, Wen-Fong and Lee, Pjek-Hwee},
{\it Maps characterized by action on zero products},
 {Pacific J. Math.}, {\bf 216}, no. 2, {(2004)}, {217--228}.

\bibitem{Din_LP:2004}
{Dinh, H. Q. and L{\'o}pez-Permouth, S. R.:}
{\it On the equivalence of codes over rings and modules},
   {Finite Fields and their Applications},
  {\bf 10}, {(4)}, {(2004)},
    {615--625}.

\bibitem{engel}
Engelking, R.: General Topology, Polish Scientific Publishers,
Warszawa (1977).

\bibitem{FerGarHer:JFS}
{Ferrer, M., Gary, M., and Hern\'andez, S.:}
{\it Representation of group isomorphisms: The compact case}
(2014), Journal of Function Spaces, Vol. 2015, Article ID 879414, 6 pages, 2015. 
doi:10.1155/2015/879414

 \bibitem{Fer_Her_Sha:Products}
  {Ferrer, M., Hern\'andez, S., and Shakhmatov, D.:}
  {\it Subgroups of direct products closely approximated by direct sums},
  (2014), Pending. 
  arXiv:1306.3954 [math.GN].

\bibitem{font-her}
Font, J.J. and Hernández, S.: \emph{On separating maps between
locally compact spaces}, Arch. Math. \textbf{63}(1994), 158-165.

\bibitem{font-her2}
Font, J.J. and Hernández, S.: \emph{Automatic continuity and
representation of certain linear isomorphisms between group
algebras}, Indag. Mathem. \textbf{6}(4)(1995), 397-409.

\bibitem{forney_trott:04} Forney, G. D. Jr. and Trott M. D.:, \textit{The
Dynamics of Group Codes: Dual Abelian Group Codes and Systems},
Proc. IEEE Workshop on Coding, System Theory and Symbolic Dynamics
(Mansfield, MA), pp. 35-65  (2004).

\bibitem{gau-jeang-wong}
Gau, H-L, Jeang, J-S, and Wong, N. C.: \emph{An algebraic approach to
the Banach-Stone theorem for separating linear bijections} Taiwanese
J. Math.  \textbf{6}(3)(2002), 399-403.

\bibitem{GluLue:2009}
{Gluesing-Luerssen, H.:}
     {\it On isometries for convolutional codes},
   {Advances in Mathematics of Communications},
    {\bf 3}, {(2)}
      {(2009)},
    {179--203}.



\bibitem{her-bec-nar}
Hernández, S., Beckenstein, E.  and Narici. L.: \emph{Banach-Stone
theorems and separating maps}, Manuscripta Math.
\textbf{86}(1995), 409-416.

\bibitem{jarosz}
Jarosz, K.: \emph{Automatic continuity of separating linear maps},
Canad. Math. Bull. \textbf{33}(2)(1990), 139-144.

\bibitem{Lau_Wong:2013}
{Lau, Anthony To-Ming and Wong, Ngai-Ching},
{\it Orthogonality and disjointness preserving linear maps between
{F}ourier and {F}ourier-{S}tieltjes algebras of locally compact groups},
{J. Funct. Anal.}, {\bf 265}, no. 4, {(2013)}, {562--593}.

\bibitem{McW:i} MacWilliams, F. J.: {\it Combinatorial problems of elementary abelian groups},
PhD thesis, Harvard University, 1962.

\bibitem{McW:ii} MacWilliams, F. J.: {\it A theorem on the distribution of weights in a systematic code},
Bell Syst. Tech. J., \textbf {42}, (1963), 79-94.

\bibitem{ohta 96}
Ohta, H.: \emph{Chains of strongly non-reflexive dual groups of
integer-valued continuous functions}, Proc. Amer. Math. Soc.
\textbf{124}(3)(1996), 961-967.

\bibitem{piret:1988}
Piret, P.: {\it Convolutional Codes; An Algebraic Approach}, MIT Press, Cambridge, MA,
1988.

\bibitem{rosenthal} Rosenthal, J., Schumacher, S. M.,  and York, E. V.: \emph{On
  Behaviors and Convolutional Codes}, IEEE Trans. on
  Information Theory, 42 (1996), no. 6, 1881--1891.

\bibitem{Stone}
Stone, M.H.: {\it Applications of the theory of boolean rings to General Topology},
Trans. Amer. Math. Soc., {\bf 41} (1937), 375-481.

\bibitem{WardWood:1996}
{Ward, H. N. and Wood, J. A.:}
     {\it Characters and the equivalence of codes},
   {Journal of Combinatorial Theory. Series A},
    {\bf 73}, no. 2, {(1996)},
    {348--352}.

\bibitem{willems:86} Willems, J. C.: \textit{From time series to linear
systems}, Parts I-III, vol. 22, pp. 561-580 and 675-694, 1986.

\bibitem{wood:2001}
Wood, J. A.: {\it The structure of linear codes
of constant weight}, Trans. Amer. Math. Soc., {\bf 354}, no. 3, (2001), 1007-1026.
\end{thebibliography}
\end{document}